\documentclass[]{interact}

\usepackage{epstopdf}

\usepackage{natbib}
\bibpunct[, ]{(}{)}{;}{a}{}{,}
\renewcommand\bibfont{\fontsize{10}{12}\selectfont}

\theoremstyle{plain}

\theoremstyle{definition}

\theoremstyle{remark}

\newtheorem{property}{Property}

\allowdisplaybreaks

\usepackage[utf8]{inputenc}

\usepackage{amsmath, amsthm, amssymb, bm,thmtools}

\usepackage{multirow}
\usepackage{algorithm}
\usepackage[noend]{algorithmic} 

\usepackage{color}

\usepackage{soul}
\usepackage{subcaption}
\usepackage{tikz-network} 
\usepackage{chemformula}
\usetikzlibrary{patterns}

\usepackage{pgfplots}
\pgfplotsset{compat=1.17}

\usepackage{pdflscape}

\def \conv{\textup{conv}}
\def \extr{\textup{extr}}

\usepackage{url}

\usepackage{bbm, amsthm}

\newcommand{\cE}{\mathcal{E}}
\newcommand{\cX}{\mathcal{X}}
\newcommand{\cK}{\mathcal{K}}
\newcommand{\cM}{\mathcal{M}}
\newcommand{\cT}{\mathcal{T}}
\newcommand{\cU}{\mathcal{U}}

\newcommand{\ii}{\mathrm{i}}

\begin{document}

\articletype{ARTICLE TEMPLATE}

\title{MIQCP and MISOCP-Based Solution Methods \\ for the Multi-Layer Thin Films Problem}

\author{
Deniz Tuncer\thanks{CONTACT Deniz Tuncer. Email: dtuncer@sabanciuniv.edu} \and
Burak Kocuk\thanks{CONTACT Burak Kocuk. Email: burak.kocuk@sabanciuniv.edu}
\affil{
\\Faculty of Engineering and Natural Sciences, Sabancı University, Istanbul, Turkey 34956
}
}

\maketitle

\begin{abstract}
The Multi-Layer Thin Films Problem is a materials science problem that aims to enhance the reflectance of a metallic substrate by designing multi-layer coatings composed of different dielectric materials and thicknesses. While previous studies on the problem mostly rely on heuristic approaches and are designed for single wavelength applications, this work addresses the problem using global optimization techniques for multiple wavelengths. We develop an exact nonconvex mixed-integer quadratically constrained programming (MIQCP) model to solve this problem. We also develop a mixed-integer second-order cone programming relaxation that has computational advantage over the MIQCP model. Our numerical experiments yield solutions that have average reflectance of 99\% over the visible spectrum (380–770 nm) and 95\% over the broad spectrum (300-3000 nm).
\end{abstract}


\maketitle

\section{Introduction}
\label{sec:Intro}

Reflectance is an electromagnetic property of materials that quantifies the percentage of incident light energy reflected by the material. In some optics applications such as telescope and laser systems, materials with high reflectance are preferred. The Multi-Layer Thin Films Problem involves a metallic surface whose reflectance is to be enhanced by coating it with thin layers of various dielectric materials.

The Multi-Layer Thin Films Problem can be considered a special case of the optimization problems of the following form \citep{kocuk2022optimization}:
\begin{subequations}\label{eq:generalproblem}
\begin{align}
\max_{T_1, \ldots, T_N,\, w} \quad & f(w) \\
\text{s.t.} \quad & p T_1 T_2 \cdots T_N = w, \label{eq:const_mult}\\
                  & T_1, \ldots, T_N \in \mathcal{T}
\end{align}
\end{subequations}

In Problem \eqref{eq:generalproblem},  \( p \in \mathbb{C}^{r \times d} \) is a given matrix that represents the initial state, 
\( f : \mathbb{C}^{r \times d} \rightarrow \mathbb{R} \) is a real-valued function, 
and \( \mathcal{T} \subseteq \mathbb{C}^{d \times d} \) is a given set of matrices. Also, \( w \in \mathbb{C}^{r \times d} \) represents the final state, after $N$ many matrices are multiplied with $p$. Since constraint \eqref{eq:const_mult} involves multiplication of matrix variables, this is a nonlinear optimization problem. This problem setting, also called switched linear systems \citep{wu2018optimal}, represents a system whose initial state is $p$, and the system evolves through each decision, $T_1,\dots, T_N$. Then, the performance of the final state $w$ of the system is evaluated by the function $f$. Despite its seemingly simple structure, even a special case of this problem, where $f$ is linear and $\cT$ is a finite set, is shown to be NP-Hard~\citep{tran2017antibiotics}. In the context of the Multi-Layer Thin Films Problem, the variables $T_1,\dots,T_N$ denote the transfer matrices corresponding to the selected coating material and thickness for layer $n$, $f$ is the reflectance function that will be introduced in Section~\ref{sec:optform} and $p$ is the $2\times2$ identity matrix.

A review of the literature reveals that the design of multi-layer thin films are typically guided by heuristic methods~\citep{larruquert99,Dobrowolski2002,hobson2004markov,shi2017optimization,keccebacs2018enhancing}. One of the earliest studies on the Multi-Layer Thin Films Problem is \cite{Turner1966}, where the authors define the problem and come up with a heuristic idea: stacking alternating layers of high and low refractive-index materials yield good reflectance values. In the literature, there are numerous studies that utilize this idea (see \cite{Southwell1980,Popov1997,keccebacs2018enhancing}). Another commonly used heuristic relies on the idea of selecting the coating thicknesses as optical thickness of a quarter wavelength (see \cite{macleod2010thin}). A special heuristic method called needle optimization is first coined by \citet{tikhonravov1996application} for this problem, and is utilized by \cite{Mashaly2024} for a practical energy saving application. In general, heuristics are utilized in the inverse design process, which involves choosing a coating, computing its reflectance and then adjusting the thicknesses (e.g., \cite{So2019}). Particle swarm optimization and numerical optimization methods are also used in a similar manner to come up with coating designs \citep{Kim2021,Zhang2024}. Recent advances in deep learning technologies have facilitated its application to the design of multi-layer thin films as well \citep{jiang2020,fouchier21}.

While heuristic methods may work well for some specific wavelengths, they achieve it with high number of layers. We suggest that optimization techniques have significant potential for achieving improved reflectance while coating with fewer layers, which may decrease the likelihood of implementation errors. To the best of our knowledge, the studies~\citet{Azunre2019,wu2018optimal,kocuk2022optimization} represent the attempts to address this problem using mathematical optimization methods. In \cite{kocuk2022optimization}, the author formulates and solves a nonconvex mixed-integer quadratically constrained program (MIQCP) for increasing the reflectance of the substrate for light coming at only a single wavelength.

In this study, our objective is to maximize the average reflectance of the substrate over a range of wavelengths. In other words, our proposed designs will be robust to changes in the wavelength of the incoming light. We aim to achieve these robust designs using global optimization methods, whose utilization in this domain is limited. Our contributions are as follows: 
We propose a novel MIQCP-based global optimization method for the Multi-Layer Thin Films Problem that aims to maximize the average reflectance over a spectrum of wavelengths. We also devise a convex mixed-integer relaxation of the MIQCP in the form of a mixed-integer second-order cone program (MISOCP). We compare our results with a heuristic from the literature and find out that our methods are able to provide better solutions that yield higher reflectance values with fewer coating layers.

Our paper is organized as follows: In Section~\ref{sec:optform}, we introduce some concepts related to the underlying optics concepts and then provide an optimization formulation. In Section~\ref{sec:solutionmethod}, we propose our MIQCP and MISOCP-based methods. In Section~\ref{sec:computation}, we present our computational results. Lastly, in Section \ref{sec:conclusion}, we discuss our key findings and the future work that may follow.

\section{Optimization Formulation}
\label{sec:optform}

Developing an optimization formulation for the Multi-Layer Thin Films Problem requires some background information on the optical properties of the materials. Therefore, before we provide a formulation, we resort to some classical optics textbooks (e.g., \cite{macleod2010thin, pedrotti2017introduction}) and introduce some basic concepts and notation used in multi-layer thin films literature.

Each material has a property called a refractive index, which is a measure of how fast the light travels through the material. This property is important because in multi-layered applications, the effect of different refractive indices determines the reflectance of the material. Let the complex refractive index of a metallic substrate such as Tungsten $s$ at a specific wavelength $\lambda$ be denoted as $\hat a_s^\lambda  \in \mathbb{C}$, where the imaginary part represents the reflection loss. Let $\cM$ represent the set of dielectric coating materials, such as Titanium Dioxide (\ch{TiO2}) and Magnesium Fluoride (\ch{MgF2}).

To quantify the reflectance of a material, we introduce the notion of transfer matrices. Transfer matrices are derived from Maxwell's equations related to the electromagnetic fields of the waves \citep{pedrotti2017introduction}. Assuming that the light is at normal incidence to the material, the transfer matrix of material $m$ with thickness $\theta$ at wavelength $\lambda$ is expressed as follows (here $\ii = \sqrt{-1}$):

\begin{equation}\label{eq:transferM}
    T_{m}^\lambda (\theta) =
     \begin{bmatrix}
    \cos  \sigma_{m,\theta}^\lambda   &  \ii \frac{\sin \sigma_{m,\theta}^\lambda}{ \hat  a_m^\lambda} \\
     \ii {\hat a_m^\lambda}{\sin \sigma_{m,\theta}^\lambda} & \cos  \sigma_{m,\theta}^\lambda
    \end{bmatrix}, \ \  \sigma_{m,\theta}^\lambda = \frac{2\pi \hat a_m^\lambda \theta}{\lambda} .
    \end{equation}

The transfer matrices of the form \eqref{eq:transferM} satisfy the following properties which will be crucial for our derivations later:

\begin{property}\label{prop:multiplicativity}
Transfer matrices are multiplicative, meaning that the effect of coating a surface with material $m_1$ of thickness $\theta_1$, and then coating it with material $m_2$ of thickness $\theta_2$ is equal to the multiplication of $T_{m_1}^\lambda (\theta_1) $ and $ T_{m_2}^\lambda (\theta_2) $.
\end{property}

\begin{property}\label{prop:determinantprop}
Determinant of transfer matrices is equal to 1, i.e., $\det(T_m^\lambda(\theta)) = 1$.
\end{property}

\begin{property}\label{prop:determinantcumulative}
Determinant of cumulative transfer matrices is equal to 1, which is a consequence of Property~\ref{prop:multiplicativity} and Property~\ref{prop:determinantprop}.
\end{property}

Notice that the transfer matrices \eqref{eq:transferM} have real values in the diagonal entries and imaginary values in the off-diagonal entries. For a given complex matrix $M \in \mathbb{C}^{2 \times 2}$, we denote an associated matrix $\tilde M  \in \mathbb{R}^{2 \times 2}$ as 
    \begin{equation*}\label{eq:tilde w def}
    \tilde M_{i,j} = \begin{cases}
    \Re(M_{ij}) & \text{ if }  (i,j)\in\{(1,1),(2,2)\} \\
    \Im(M_{ij}) & \text{ if }  (i,j)\in\{(1,2),(2,1)\} 
    \end{cases}.
    \end{equation*}

Since the cumulative effect of applying multiple thin films is governed by multiplication of their associated transfer matrices, the resulting multiplication is called a cumulative transfer matrix. Now, suppose that we have a multi-layer thin film with the cumulative transfer matrix $w\in\mathbb{C}^{2\times2}$ coated on substrate $s$. Then, the reflectance of the coated material at  wavelength $\lambda$ is computed as 
\begin{equation}\label{eq:reflectance}
 R_s^\lambda (w) := 
 \frac{ (\tilde w_{11}-\Im({\hat a_s^\lambda}) \tilde w_{12}-\Re({\hat a_s^\lambda})\tilde w_{22} )^2 + (\tilde w_{21}+\Im({\hat a_s^\lambda})\tilde w_{22}-\Re({\hat a_s^\lambda) \tilde w_{12}})^2 }
 { (\tilde w_{11}-\Im({\hat a_s^\lambda})\tilde w_{12}+\Re({\hat a_s^\lambda})\tilde w_{22} )^2 + (\tilde w_{21}+\Im({\hat a_s^\lambda})\tilde w_{22}+\Re({\hat a_s^\lambda) \tilde w_{12}})^2 }.
\end{equation}

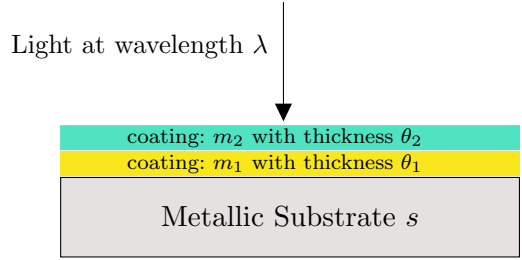
\begin{figure}[H]
    \centering
    \begin{tikzpicture}
[x=0.75pt,y=0.75pt,yscale=-1,xscale=1]
\tikzset{every picture/.style={line width=0.75pt}} 
\draw  [fill={rgb, 255:red, 230; green, 225; blue, 225 }  ,fill opacity=1 ] (215,115) -- (444,115) -- (444,155) -- (215,155) -- cycle ;
\draw  [color={rgb, 255:red, 248; green, 231; blue, 28 }  ,draw opacity=1 ][fill={rgb, 255:red, 248; green, 231; blue, 28 }  ,fill opacity=1 ] (215,102) -- (444,102) -- (444,114) -- (215,114) -- cycle ;
\draw  [color={rgb, 255:red, 80; green, 227; blue, 194 }  ,draw opacity=1 ][fill={rgb, 255:red, 80; green, 227; blue, 194 }  ,fill opacity=1 ] (215,89) -- (444,89) -- (444,101) -- (215,101) -- cycle ;
\draw    (326.2,26.8) -- (326.2,47.8) -- (326.2,83.8) ;
\draw [shift={(326.2,86.8)}, rotate = 270] [fill={rgb, 255:red, 0; green, 0; blue, 0 }  ][line width=0.08]  [draw opacity=0] (8.93,-4.29) -- (0,0) -- (8.93,4.29) -- cycle    ;

\draw (329.61,134.5) node   [align=left] {Metallic Substrate $s$};
\draw (247,102) node [anchor=north west][inner sep=0.75pt]  [font=\footnotesize] [align=left] {coating: $\displaystyle m_{1}$ with thickness $\displaystyle \theta _{1}$ };
\draw (247,89) node [anchor=north west][inner sep=0.75pt]  [font=\footnotesize] [align=left] {coating: $\displaystyle m_{2}$ with thickness $\displaystyle \theta _{2}$ };
\draw (189,42) node [anchor=north west][inner sep=0.75pt]  [font=\small] [align=left] {Light at wavelength $\displaystyle \lambda $};
\end{tikzpicture}
    \caption{Illustration of two layers of coating on top of a substrate.}
    \label{fig:thinfilms}
\end{figure}

To make the problem and the concepts more concrete, we provide an illustration of a two-layer thin film in Figure~\ref{fig:thinfilms}.

Assume that on top of the metallic substrate, we coat a thin film of material $m_1$ with thickness $\theta_1$ as the first layer. As the second layer, we add a thin film of material $m_2$ with thickness $\theta_2$. Then, light beams at wavelength $\lambda$ are directed onto the surface and the reflectance is measured. Due to the multiplicative property, we can calculate the cumulative transfer matrix ($w$) in wavelength $\lambda$ by multiplying the associated transfer matrices $T_{m_1}^\lambda(\theta_1)$ and $T_{m_2}^\lambda(\theta_2)$. Then, the cumulative transfer matrix can be used to calculate the reflectance of the material as in \eqref{eq:reflectance}.

\subsection{Generic Optimization Formulation}\label{sec:genericopt}

Suppose that we aim to coat a substrate $s$ with $N$ many thin film layers so that its reflectance is maximized for a finite set of wavelengths denoted by $\Lambda$. We would like to design a multi-layer coating that maximizes the average reflectance when the light comes at any wavelength $\lambda \in \Lambda$. Let us denote the intensity of light at wavelength $\lambda$ as $\phi^\lambda$, $\lambda\in\Lambda$. Let $\cM$ be the set of dielectric coating materials and $\Theta_m$ be the set of admissible thickness values for material $m\in\cM$. 

In order to formulate this situation as an optimization problem, we define decision variables $\mu_n$ and $\theta_n$ representing the material and its thickness selected for layer $n$, $n=1,\dots,N$, and present the following generic formulation:
\begin{equation} \label{eq:generic formulation}
\max_{\mu, \theta} \left\{  \sum_{\lambda \in \Lambda} \phi^\lambda R_s^\lambda \left( \prod_{n=1}^N T_{\mu_n}^\lambda (\theta_n) \right)  :
  \mu_n \in \cM ,  \   \theta_n \in \Theta_{\mu_n}   , \   n=1,\dots,N  \right\} .
\end{equation}
Note that formulation~\eqref{eq:generic formulation} is nonlinear due to the product of matrices and the definition of the reflectance function \eqref{eq:reflectance}, and also contains discrete decisions.

We will adapt some ideas from~\citet{kocuk2022optimization} to ``bilinearize" and ``linearize" the above formulation. Let us first define variable matrices $M_n^\lambda \in \mathbb{C}^{2\times2}$ and $u_n^\lambda \in \mathbb{C}^{2\times2}$  to represent the transfer matrix of layer $n$ at wavelength $\lambda$ and the cumulative transfer matrix up to layer $n$ at wavelength $\lambda$, respectively. 
Due to Proposition~3.2 in \citet{kocuk2022optimization}, $ R_s^\lambda ( w)$ defined in~\eqref{eq:reflectance} can be alternatively represented as
 $
 R_s^\lambda ( w) = 1-\frac{ 4 \Re({\hat a_s^\lambda}) } { D_s^\lambda ( w)  }
,
$
where
\begin{equation}\label{eq:defDen}
D_s^\lambda ( w) := (\tilde w_{11}-\Im({\hat a_s^\lambda})\tilde w_{12})^2 +  (\Re({\hat a_s^\lambda) \tilde w_{12})^2 +
  (\tilde w_{21}+\Im({\hat a_s^\lambda})\tilde w_{22}})^2 
  + (\Re({\hat a_s^\lambda})\tilde w_{22})^2  + 2\Re({\hat a_s^\lambda}) .
  \end{equation} 
  After defining $u_0^\lambda:=I$ and $w^\lambda:=u_N^\lambda$ for convenience,  we  obtain the following formulation that is equivalent to \eqref{eq:generic formulation}:  
\begin{subequations} \label{eq:generic reformulation}
\begin{align}
\max_{\mu,\theta, M, u, w} \ &  \sum_{\lambda \in \Lambda}  \phi^\lambda \left(1 - \frac{ 4\Re({\hat a_s^\lambda})}{D_s^\lambda \big( w^\lambda \big)  } \right)\label{eq:generic reform obj} \\
\text{s.t.} 
\ &  M_n^\lambda = T_{\mu_n}^\lambda (\theta_n)   &  \lambda&\in\Lambda;  \mu_n \in \cM ,    \theta_n \in \Theta_{\mu_n}   ,    n=1,\dots,N  \label{eq:generic reform matrix} \\
\ & u_0^\lambda = I,  u_N^\lambda = w^\lambda &  \lambda&\in\Lambda \label{eq:initial-final}   \\
\ & u_{n-1}^\lambda M_n^\lambda = u_n^\lambda &  \lambda&\in\Lambda, n=1,\dots,N . \label{eq:recursion}   
\end{align}
\end{subequations}

We will specialize this formulation for the case when the applicable thicknesses for coatings are discrete sets ($\Theta_m$ is a finite set for $m\in\cM$) and we consider maximizing the average reflectance over multiple wavelengths (i.e., $|\Lambda|>1$).

\subsection{Discrete Material Thicknesses}
\label{sec:form-disc}

Assume that the thickness of each layer can be selected from a finite set $\Theta_m$ for each $m \in \cM$. Suppose that we compute \textit{fixed} transfer matrices $\hat T_{m,\theta}^\lambda := T_m^\lambda (\theta)$ a priori for each $\lambda\in\Lambda$, $m\in\cM$ and $\theta\in\Theta_m$. We define a set of binary decision variables $x_{n,(m,\theta)}$, which take value one if layer $n$ is coated with material $m$ of thickness $\theta$, and zero otherwise. 
Further suppose that we have polyhedral sets $\cU^\lambda_n$, which contain all possible values of the decision variable $u_n^\lambda$, at our disposal. These polyhedral sets can be constructed using the recursive relation \eqref{eq:recursion}, and they are  bounded since the number of transfer matrices are finite. In order to linearize the relation~\eqref{eq:recursion}, we use the disjunctive formulation proposed in~\cite{kocuk2022optimization} by exploiting the finiteness of transfer matrices. Let us define auxiliary variables $v_{n-1, (m,\theta)}^\lambda$, which take value $u_{n-1}^\lambda$ if the binary variable $x_{n,(m,\theta)}$ takes value one and zero otherwise. Then, we obtain the following equivalent optimization problem:

    \begin{subequations} \label{eq:discrete}
    \begin{align}
    \max_{\mu,\theta, M, u, w} \ &  \sum_{\lambda \in \Lambda}  \phi^\lambda \left(1 - \frac{ 4\Re({\hat a_s^\lambda})}{D_s^\lambda \big( w^\lambda \big)  } \right)\label{eq:generic reform obj disc} \\
    \text{s.t.} \ &  \eqref{eq:initial-final} \notag \\
    \ & \sum_{m\in\cM} \sum_{\theta\in\Theta_m} v^\lambda_{n-1,(m,\theta)} = u^\lambda_{n-1}  & n&=1,\dots,N \label{eq:disj1}  \\
    \ &  v^\lambda_{n-1,(m,\theta)} \hat T_{m,\theta}^\lambda  = u^\lambda_{n}  & n&=1,\dots,N \label{eq:disj2}  \\
    \ &  v^\lambda_{n-1,(m,\theta)} \in \cU_{n-1}^\lambda x_{n,(m,\theta)}    & n&=1,\dots,N,  \lambda\in\Lambda , m\in\cM, \theta\in\Theta_m \label{eq:disj3}  \\
    \ & \sum_{m\in\cM} \sum_{\theta\in\Theta_m} x_{n,(m,\theta)} = 1  & n&=1,\dots,N \label{eq:one material-thickness}  \\
    \ & x_{n,(m,\theta)} \in \{0,1\} & n&=1,\dots,N , m\in\cM, \theta\in\Theta_m.  \label{eq:binary-x}
    \end{align}
    \end{subequations}

We can further reformulate the objective function \eqref{eq:generic reform obj disc} using two more sets of auxiliary variables $d_s^\lambda$ and $f_s^\lambda$ satisfying the following:
\begin{subequations} \label{eq:define d and f}
\begin{align}
  & \  f_s^\lambda d_s^\lambda \ge 4 {\Re({\hat a_s^\lambda})}  &\lambda& \in\Lambda  \label{eq:define d and f-SOC} \\ 
  & \  d_s^\lambda \le D_s^\lambda \big( w^\lambda \big)   &\lambda& \in\Lambda \label{eq:define d and f-reverse}\\ 
  & \  f_s^\lambda \ge  0 , d_s^\lambda \ge  0    &\lambda& \in\Lambda.\label{eq:signs for d and f}
\end{align}
\end{subequations}
We observe that inequality~\eqref{eq:define d and f-SOC} is second-order cone representable, which is an important fact since it will help us build a convex relaxation of the problem later. Also, inequality~\eqref{eq:define d and f-reverse} is reverse convex and is the only remaining nonconvex relation (except the binary restrictions). 
The resulting nonconvex MIQCP is defined as below:
\begin{equation}\label{eq:discrete-nonconvex}
    \max_{u, v, w, x, f,d} \left\{  \sum_{\lambda \in \Lambda}  \phi^\lambda \left(1 - {f_s^\lambda  }\right)  :  \eqref{eq:initial-final}  , \eqref{eq:disj1} - \eqref{eq:binary-x} , \eqref{eq:define d and f} \right\}.
\end{equation}

\section{Solution Method}
\label{sec:solutionmethod}

Two reasonable approaches appear to exist for solving the thin film spectrum optimization problem with discrete thicknesses:

\begin{enumerate}
    \item MIQCP-based Approach: We solve the nonconvex MIQCP problem~\eqref{eq:discrete-nonconvex} directly using an off-the-shelf global optimization solver such as Gurobi \citep{gurobi}.
    \item MISOCP-based Approach: 
    We first find linear overapproximators of $D_s^\lambda (w^\lambda) $ over the variable bounds obtained from Algorithm~\ref{alg:boundtightening} and  $\det(w^\lambda)=1$ (recall that cumulative transfer matrices have their determinant equal to one). Using these linear overapproximators, we   build an MISOCP relaxation of the original nonconvex MIQCP problem as problem~\eqref{eq:misocp_problem}, which has a computational advantage due to its convex structure. 
\end{enumerate}

In the remainder of this section, we propose a bound tightening procedure that is utilized in both MIQCP and MISOCP-based methods in Section~\ref{sec:bound} and then develop the MISOCP formulation in Section~\ref{sec:approximation}.

\subsection{Bound Tightening}
\label{sec:bound}

The effectiveness of global optimization methods relies heavily on the presence of tight variable bounds, which also effect the quality of the solutions obtained by the solvers. Let us denote the lower bound and the upper bound for the $(i,j)$'th entry of the variable matrices ($u_n^{\lambda},n=1,\dots,N$), as $\underline{u}_{n,(i,j)}^\lambda$ and $\overline{u}_{n,(i,j)}^\lambda$, respectively. 
Then, we are able to present our bound tightening algorithm in Algorithm~\ref{alg:boundtightening}.

\begin{algorithm}[H]
\caption{Bound Tightening Algorithm.}
\label{alg:boundtightening}
\begin{algorithmic}[1]
\STATE Set $n = 1$.
\FORALL{$\lambda \in \Lambda$}
\STATE Define the set of bound matrices, initializing $\mathcal{B}_{0}^{\lambda} = \{I\}$.
\WHILE{$n \leq N$}
\STATE $\mathcal{C}_{n}^{\lambda} = \{~\} $
\FORALL{$\hat{T}^\lambda_{m,\theta}, m\in \mathcal{M}, \theta \in \Theta_m$ and $B \in \mathcal{B}_{n-1}^{\lambda} $}
\STATE $\mathcal{C}_{n}^{\lambda} = \mathcal{C}_{n}^{\lambda} \cup \{B \hat{T}^\lambda_{m,\theta}\}$
\ENDFOR
\FORALL {$i\in \{1,2\}$, $j \in \{1,2\}$} \STATE Set  $\underline{u}_{n,(i,j)}^\lambda =$ $\min_{C \in \mathcal{C}_{n}^{\lambda}} C_{n,(ij)}^{\lambda}$ and $\overline{u}_{n,(i,j)}^\lambda =$ $\max_{C \in \mathcal{C}_{n}^{\lambda}} C_{n,(ij)}^{\lambda}$ .\label{lst:line:boundsetting}
\ENDFOR
\STATE Set $\mathcal{B}_n^\lambda:=\left\{\begin{bmatrix}
u_{11}& u_{12} \\
u_{21} & u_{22}
\end{bmatrix} : u_{ij} \in \{\underline{u}_{n,(i,j)}^\lambda,\overline{u}_{n,(i,j)}^\lambda \} , i\in\{1,2\},j\in\{1,2\}\right\}$. 
\STATE $n=n+1$.
\ENDWHILE
\ENDFOR
\end{algorithmic}
\end{algorithm}

Note that variables $u_n^\lambda$ are complex numbers; we apply $\min$ and $\max$ operators separately for the real and imaginary parts of the variables in  line~\ref{lst:line:boundsetting} of Algorithm~\ref{alg:boundtightening}.

Recall that $u_N^\lambda = w^{\lambda} $ for all $\lambda \in \Lambda$. Therefore, variable bounds found after running Algorithm~\ref{alg:boundtightening} are valid for $w^\lambda$.

\subsection{Relaxation by Overapproximation of $D_s^\lambda (w^\lambda) $ }
\label{sec:approximation}

Our aim is to find a linear overapproximator of the convex quadratic function $D_s^\lambda (w^\lambda) $ over the set defined by variable bounds and the constraint $\det(w^\lambda)=1$.  In order to simplify the notation, we will denote variables $\tilde w_{11}^\lambda, \tilde w_{22}^\lambda, \tilde w_{12}^\lambda, \tilde w_{21}^\lambda$ by $x_1, x_2, x_3, x_4$ and function $D_s^\lambda (w^\lambda)$ by $g(x)$.

Given the variable bounds, $\underline x, \overline x \in \mathbb{R}^4$ with $\underline x \le \overline x$ and a convex function $g:\mathbb{R}^4 \to \mathbb{R}$, our aim is to find a linear overapproximator $\alpha_0 + \sum_{\ell=1}^4 \alpha_\ell x_\ell$ of $g(x)$ over $ \cX := \{  x\in [\underline x , \overline x] :  x_{1} x_{2} + x_{3} x_{4} = 1\}$. In other words, we aim to find $\alpha \in \mathbb{R}^5$ such that
\[
\alpha_0 + \sum_{\ell=1}^4  \alpha_{\ell} x_{\ell} \ge g(x) \quad \text{ for all } x \in \cX.
\]
Since $g$ is convex, it suffices to satisfy the above condition at the extreme points of set $\cX$, which are characterized as follows \citep{kocuk2018matrix}:  
\begin{equation*}
\cE :=  \left(  \bigcup_{\substack{\hat x_1\in\{\underline x_1,\overline x_1\}  \\ \hat x_2\in\{\underline x_2,\overline x_2\}  } }  (\hat x_1, \hat x_2) \times \extr \big( \cX^{34} ( \hat x_{1}, \hat x_{2}) \big ) \right)
    \bigcup
  \left (  \bigcup_{ \substack{\hat x_3\in\{\underline x_3,\overline x_3\} \\ \hat x_4\in\{\underline x_4,\overline x_4\} } }  \extr \big( \cX^{12} ( \hat x_{3}, \hat x_{4}) \big ) \times (\hat x_3, \hat x_4) \right ),
\end{equation*}
where
  $$  \cX^{12} ( \hat x_{3}, \hat x_{4})  := \{ (x_{1}, x_{2}) \in [\underline x_1, \overline x_1] \times [\underline x_2, \overline x_2]:     x_{1} x_{2} = 1 - \hat x_{3} \hat x_{4}    \} , $$ and
  $$  \cX^{34} ( \hat x_{1}, \hat x_{2})  := \{ (x_{3}, x_{4}) \in [\underline x_3, \overline x_3] \times [\underline x_4, \overline x_4]:     x_{3} x_{4} = 1 - \hat x_{1} \hat x_{2}    \}  . $$ 

Now, we study the characteristics of sets $\cX^{12} ( \hat x_{3}, \hat x_{4})$ and $\cX^{34} ( \hat x_{1}, \hat x_{2})$. In general, the convex hull of a nonempty set $ K := \{ (y_1, y_2) \in [\underline y_1, \overline y_1] \times [\underline y_2, \overline y_2]:     y_{1} y_{2} = \beta   \} $ is one of the following types:
\begin{itemize}
    \item Polytope: This is when variable bound constraints and constraint $y_1 y_2 = \beta$ intersects both at the negative and positive orthant. In this case, extreme points of the set $\conv (K)$ consists of at most four elements given by $  \{(\underline{y}_1,\beta/\underline{y}_1),\allowbreak (\overline{y}_1, \beta/\overline{y}_1), \allowbreak(\beta/\underline{y}_2, \underline{y}_2) ,\allowbreak (\beta/\overline{y}_2,\overline{y}_2) \}$.

    For example, assume that the set $K$ is characterized by $\beta=4$ and $\underline y_1 =-3, ~\overline y_1=3,~ \underline y_2=-2,~ \overline y_2=2$. Then, the extreme points of set $\conv(K)$ (solid dots) and the set $\conv(K)$ (shaded region) can be shown as in Figure~\ref{fig:poly-rep-case}.

    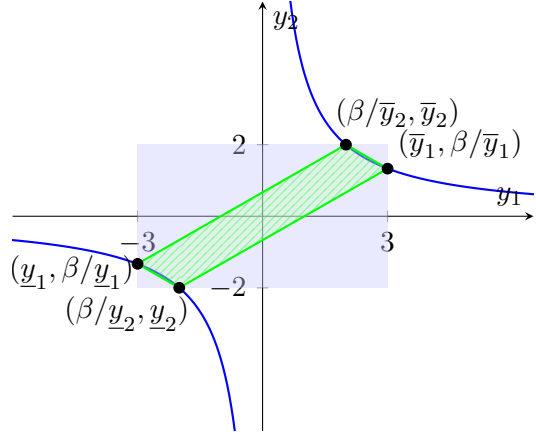
\begin{figure}[H]
    \centering
\begin{tikzpicture}
  \begin{axis}[
    axis lines=middle,
    xlabel={$y_1$},
    ylabel={$y_2$},
    xmin=-6, xmax=6.5,
    ymin=-6, ymax=6,
    samples=200,
    domain=-10:10,
    restrict y to domain=-10:10,
    grid=none,
    xtick={-3,0,3},
    ytick={-2,0,2},
  ]

    \addplot [
      draw=none,
      fill=blue!20,
      fill opacity=0.4,
    ] coordinates {
      (-3,-2) (3,-2) (3,2) (-3,2) 
    } -- cycle;

    \addplot[blue, thick, domain=-10:-0.05] {4/x};
    \addplot[blue, thick, domain=0.05:10] {4/x};

    \addplot[mark=*] coordinates {(2,2)};
    \node[anchor=south west] at (axis cs:1.5,2.2) {$(\beta/\overline{y}_2, \overline{y}_2)$};

    \addplot[mark=*] coordinates {(3,1.33)};
    \node[anchor=south west] at (axis cs:3,1.33) {$(\overline{y}_1, \beta/\overline{y}_1)$};

    \addplot[mark=*] coordinates {(-2,-2)};
    \node[anchor=south west] at (axis cs:-5,-3.5) {$(\beta/\underline{y}_2, \underline{y}_2)$};

    \addplot[mark=*] coordinates {(-3,-1.33)};
    \node[anchor=south west] at (axis cs:-6.33,-2.33) {$(\underline{y}_1, \beta/\underline{y}_1)$};

\addplot[
  draw=green,
  thick,
  fill=green!20,
  fill opacity=0.4
] coordinates {
  (-3,-1.33)
  (2,2)
  (3,1.33)
  (-2,-2)
} -- cycle;

\addplot[
    draw=green,
    pattern color=green,
    pattern=north east lines,
    opacity=0.6
    ] coordinates {
  (-3,-1.33)
  (2,2)
  (3,1.33)
  (-2,-2)
};

  \end{axis}
\end{tikzpicture}

    \caption{An example case when $\conv(K)$ is polyhedrally representable.}
    \label{fig:poly-rep-case}
\end{figure}
    \item Second-order cone representable:
   This is when variable bound constraints and $y_1 y_2 = \beta$ constraint intersect only at negative or the positive orthant. This results in a second-order cone representable set $\conv (K)$. Hence, there are infinitely many extreme points in this case. However, for practical purposes, we can outer-approximate this set with a polytope and obtain its finite list of extreme points instead. For example, assume that the set $K$ is characterized by $\beta=4$ and $\underline y_1 =-1, ~\overline y_1=3,~ \underline y_2=-2,~ \overline y_2=4$. Then, the points where the variable bound constraints and the constraint  $y_1 y_2 = \beta$ intersect are as shown (solid dots) in Figure~\ref{fig:soc-rep-case}.
   
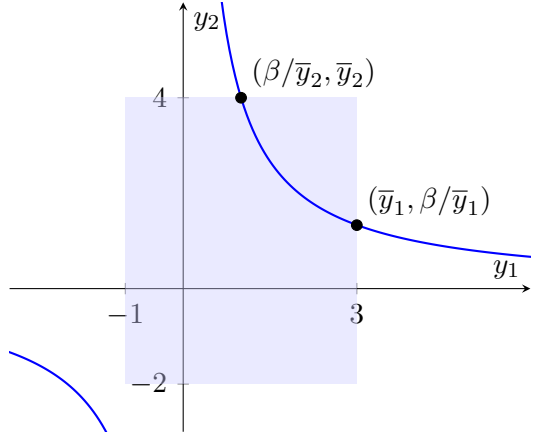
\begin{figure}[H]
    \centering
\begin{tikzpicture}
  \begin{axis}[
    axis lines=middle,
    xlabel={$y_1$},
    ylabel={$y_2$},
    xmin=-3, xmax=6,
    ymin=-3, ymax=6,
    samples=200,
    domain=-10:10,
    restrict y to domain=-10:10,
    grid=none,
    xtick={-1,0,3},
    ytick={-2,0,4},
  ]

    \addplot [
      draw=none,
      fill=blue!20,
      fill opacity=0.4,
    ] coordinates {
      (-1,-2) (3,-2) (3,4) (-1,4) 
    } -- cycle;

    \addplot[blue, thick, domain=-10:-0.05] {4/x};
    \addplot[blue, thick, domain=0.05:10] {4/x};

    \addplot[mark=*] coordinates {(1,4)};
    \node[anchor=south west] at (axis cs:1,4) {$(\beta/\overline{y}_2, \overline{y}_2)$};

    \addplot[mark=*] coordinates {(3,1.33)};
    \node[anchor=south west] at (axis cs:3,1.33) {$(\overline{y}_1, \beta/\overline{y}_1)$};

  \end{axis}
\end{tikzpicture}
    \caption{An example case when $\conv(K)$ is second-order cone representable.}
    \label{fig:soc-rep-case}
\end{figure}

Assume that the variable bound constraints and $y_1 y_2=\beta$ constraint intersect at the positive orthant, i.e. $\overline{y}_1 \overline{y}_2 \geq \beta$. This is without loss of generality, since for the cases when the intersection is another orthant, we can handle it by negating the necessary variable bounds. Then, we can update the bounds we will use to outerapproximate this set as follows: 
    \begin{itemize}
        \item $\underline{y}_1 = \min \{\overline{y}_1,\beta/\overline{y}_2 \}$, $\overline{y}_1 = \max \{ \overline{y}_1,\beta/\overline{y}_2 \}$,
        \item $\underline{y}_2 = \min \{ \beta/\overline{y}_1, /\overline{y}_2 \}$, $\overline{y}_2 = \max \{ \beta/\overline{y}_1, /\overline{y}_2 \}$.
    \end{itemize}
    
Set $K$ consists of two points that satisfy the constraint $y_1y_2 = \beta$ at the boundaries of the box defined by the variable bounds. To find a polyhedral outerapproximating set of $K$, we first calculate the geometric mean of the bounds of the first variable (denoted with filled square in Figure~\ref{fig:soc-rep-case-yellowbox}) as $(y_1,y_2) = (\sqrt{|\underline{y}_1 \overline{y}_1|}, \beta/ \sqrt{|\underline{y}_1 \overline{y}_1|})$. We then draw the tangent line to the constraint $y_1y_2 = \beta$ at that point and find the intersection points of the tangent line with the new variable bounds. Then, an outerapproximating set of $\conv(K)$ is the shaded area in Figure~\ref{fig:soc-rep-case-yellowbox}.

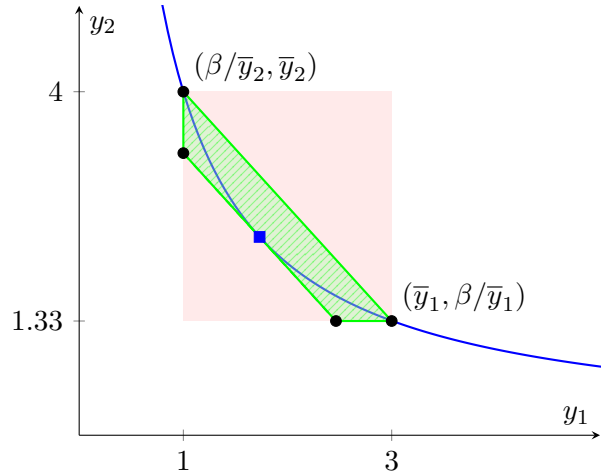
\begin{figure}[H]
\centering
\begin{tikzpicture}
  \begin{axis}[
    axis lines=middle,
    xlabel={$y_1$},
    ylabel={$y_2$},
    xmin=0, xmax=5,
    ymin=0, ymax=5,
    samples=200,
    domain=-10:10,
    restrict y to domain=-10:10,
    grid=none,
    xtick={1,3},
    ytick={1.33,4},
  ]

    \addplot [
      draw=none,
      fill=red!20,
      fill opacity=0.4,
    ] coordinates {
      (1,1.33) (3,1.33) (3,4)  (1,4)
    } -- cycle;

    \addplot[blue, thick, domain=-10:-0.05] {4/x};
    \addplot[blue, thick, domain=0.05:10] {4/x};

    \addplot[mark=*] coordinates {(1,4)};
    \node[anchor=south west] at (axis cs:1,4) {$(\beta/\overline{y}_2, \overline{y}_2)$};

    \addplot[mark=*] coordinates {(3,1.33)};
    \node[anchor=south west] at (axis cs:3,1.33) {$(\overline{y}_1, \beta/\overline{y}_1)$};

    \addplot[mark=*] coordinates {(1,3.285)};

    \addplot[mark=*] coordinates {(2.465,1.33)};

    \addplot[mark=square*,color=blue,] coordinates {(1.732,2.309)};
    
    \addplot[
    dashed,
    thick,
    color=green
] coordinates {
    (1,3.285)
    (2.465,1.33)
};

\addplot[
    draw=green,
    thick,
    fill=green!30,
    fill opacity=0.4
] coordinates {
    (1,4)
    (3,1.33)
    (2.465,1.33)
    (1,3.285)
} -- cycle;

\addplot[
    draw=green,
    pattern color=green,
    pattern=north east lines,
    opacity=0.6
    ] coordinates {
    (1,4)
    (3,1.33)
    (2.465,1.33)
    (1,3.285)
};

  \end{axis}
\end{tikzpicture}
    \caption{Outerapproximation when $\conv(K)$ second-order cone representable.}
    \label{fig:soc-rep-case-yellowbox}
\end{figure}

\end{itemize}

Let us denote the set of all points that are collected as $\cK \subseteq \mathbb{R}^4$. Using these points, we aim to obtain overapproximating functions whose coefficients $\alpha$ satisfy $\alpha_0 + \sum_{\ell=1}^4 \alpha_\ell x_\ell \geq g(x)$. Let us denote the set of all 5-element subsets of $\cK$ as $\cK^5:= \{ \cT \subseteq \cK :|\cT|=5 \}$. For each subset $\mathcal{T} = \{k^{(1)}, \dots, k^{(5)}\} \in \mathcal{K}^{5}$,
we determine the coefficients $\alpha$ of the affine function passing through these points
by solving the following linear system:

\begin{equation}\label{eq:linearsystem}
\begin{bmatrix}
    1 & k_1^{(1)} & k_2^{(1)} & k_3^{(1)}  & k_4^{(1)} \\
    1 & k_1^{(2)} & k_2^{(2)} & k_3^{(2)}  & k_4^{(2)} \\
    \vdots & \vdots & \vdots & \vdots &\vdots\\
    1 & k_1^{(5)} & k_2^{(5)} & k_3^{(5)}  & k_4^{(5)}
\end{bmatrix} \begin{bmatrix} \alpha_0 \\ \alpha_1 \\ \alpha_2 \\ \alpha_3 \\ \alpha_4\end{bmatrix} = \begin{bmatrix} g(k^{(1)}) \\ g(k^{(2)}) \\  g(k^{(3)}) \\ g(k^{(4)}) \\ g(k^{(5)})\end{bmatrix}.
\end{equation}
Then, we can find overapproximators for $D_s^\lambda(w^\lambda)$ as explained in Algorithm~\ref{alg:hyperplanefinder} for each wavelength $\lambda$ separately to construct  $\mathcal{H}_\lambda$ sets. 
\begin{algorithm}[H]
\caption{}
\label{alg:hyperplanefinder}
\begin{algorithmic}
\STATE Initialize $\mathcal{H}_\lambda = \emptyset$ (set of overapproximating functions).
\FORALL{$\cT \in \cK^5$}
\STATE Try to solve the linear system \eqref{eq:linearsystem}.
\STATE Obtain the coefficient vector $\alpha \in \mathbb{R}^5$.
\IF{$\alpha_0 + \sum_{\ell=1}^4 \alpha_\ell k_\ell \geq g(k)$ for all $k \in \cK$}
        \STATE Add the function defined by $\alpha$ to $\mathcal{H}_\lambda$.
\ENDIF
\ENDFOR
\end{algorithmic}
\end{algorithm}

Then, we can replace constraint \eqref{eq:define d and f-reverse} by the following constraint:

\begin{equation}
\label{eq:overapprox}
\alpha_0^{(h)} + \alpha_1^{(h)} \tilde{w}^\lambda_{11}+ \alpha_2^{(h)} \tilde{w}^\lambda_{22} + \alpha_3^{(h)} \tilde{w}^\lambda_{12}+ \alpha_4^{(h)} \tilde{w}^\lambda_{21} \geq d_s^\lambda 
\qquad
 \lambda \in \Lambda,\;
 h \in \mathcal{H}_\lambda.
\end{equation}

Then, we can write the MISOCP model as follows:

\begin{equation}\label{eq:misocp_problem}
    \max_{u, v, w, x, f,d} \left\{  \sum_{\lambda \in \Lambda}  \phi^\lambda \left(1 - {f_s^\lambda  }\right)  :  \eqref{eq:initial-final}  , \eqref{eq:disj1} - \eqref{eq:binary-x} , \eqref{eq:define d and f-SOC},\eqref{eq:signs for d and f}, \eqref{eq:overapprox} \right\}.
\end{equation}

\section{Computational Experiments}\label{sec:computation}

\subsection{Computational Setting}\label{subsec:compsetting}

We provide the results of our computational experiments in this section. We use a 64-bit workstation with two Intel(R) Xeon(R) Gold 6248R CPU (3.00GHz) processors (256 GB RAM) and the Python programming language. We utilize Gurobi 11 to solve the nonconvex MIQCPs. We set the parameter \texttt{BestObjStop} as 0.995, \texttt{MIPGap} as 0.01 and \texttt{TimeLimit} as 18000 seconds (unless otherwise stated). As the optimality gap improves very slowly over time, we focus on finding good quality feasible solutions and set the \texttt{MIPFocus} parameter as 1 and store up to 10 best solutions in the solution pool. If the solver finds an optimal solution within the time limit, it utilizes the remaining time for investigating other good quality solutions to populate the solution pool.

Following the observations from \cite{kocuk2022optimization}, we set the coating materials as $\mathcal{M} = \{\ch{TiO2},\ch{MgF2} \}$. We consider Molybdenum, Niobium, Tantalum and Tungsten substrates in our study.

We carry out some preliminary experiments and observe that Titanium Dioxide (\ch{TiO2}) is used for the odd numbered layers and Magnesium Fluoride (\ch{MgF2}) for the even numbered layers in all optimal solutions. To achieve tighter variable bounds, we incorporate this information into Algorithm~\ref{alg:boundtightening} by restricting the set of coating materials based on the layer index: specifically, $\mathcal{M}^n =\{$\ch{TiO2}$\}$ if $n$ is odd and $\mathcal{M}^n =\{$\ch{MgF2}$\}$ if $n$ is even.

The refractive indices that are required as an input to calculate the transfer matrices of the material-thickness pairs are obtained from \cite{keccebacs2018enhancing}, taken from multiple sources \citep{malitson1965interspecimen,palik1998handbook,dodge1984refractive,dodge2refractive,golovashkin1969optical}.

\begin{figure}[H]
\centering
    \begin{subfigure}{0.45\linewidth}
\includegraphics[width=\linewidth]{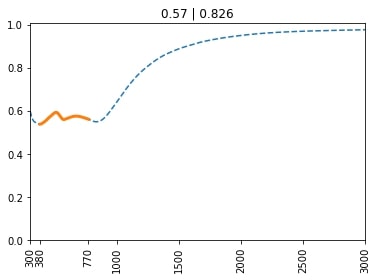}
\caption{Molybdenum}
\label{fig:subfig1bnat}
    \end{subfigure}\hfil
    \begin{subfigure}{0.45\linewidth}
\includegraphics[width=\linewidth]{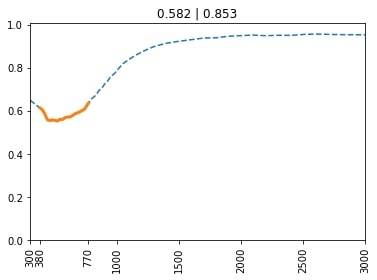}
\caption{Niobium}
\label{fig:subfig2bnat}
    \end{subfigure}

\smallskip
     \begin{subfigure}{0.45\linewidth}
\includegraphics[width=\linewidth]{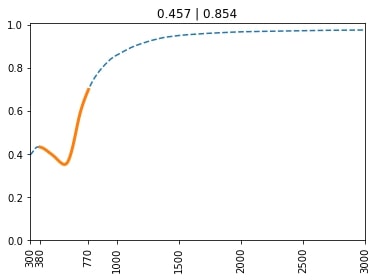}
\caption{Tantalum}
\label{fig:subfig3bnat}
    \end{subfigure}\hfil
     \begin{subfigure}{0.45\linewidth}
\includegraphics[width=\linewidth]{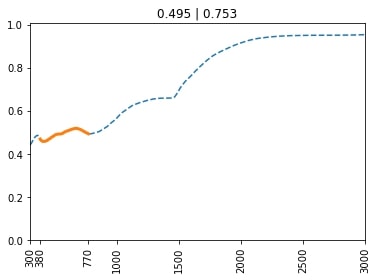}
\caption{Tungsten}
\label{fig:subfig4bnat}
    \end{subfigure}
    \caption{Reflectance profiles of uncoated substrates.}
\label{fig:natural_reflectance_profiles}
\end{figure}

We present the reflectance profiles of the uncoated substrates in Figure~\ref{fig:natural_reflectance_profiles}, where the solid line represents the reflectance values in the visible spectrum (380-770 nm) while the dashed line represents the reflectance values in the broad spectrum (300-3000 nm).

The subfigure titles in Figures~\ref{fig:natural_reflectance_profiles} and \ref{fig:vis_reflectance_profiles} present the average reflectance values calculated over the visible spectrum and the broad spectrum, respectively. For example, for the Molybdenum substrate (Figure~\ref{fig:subfig1bnat}), the average reflectance for the visible spectrum and the broad spectrum are 0.570 and 0.826, respectively. 

We would like to emphasize that the reported reflectance values do not correspond to the optimal objective function values. The optimization models are first solved to obtain an optimal solution, which is then evaluated on a finer grid to assess its performance and the resulting values are reported.

\subsection{Experiments for Single Wavelength}\label{sec:singlewavelength}

In this section, we solve the problem with MIQCP and MISOCP-based methods for some wavelengths in the visible spectrum with six layers of coating. In the experiments, the allowed thicknesses for coating are $\Theta_{\ch{TiO2}} := \{20,30,\dots,140\}$ and $\Theta_{\ch{MgF2}} := \{50,60,\dots,280\}$. Resulting reflectance values and running times in seconds are presented in Table \ref{tab:compare_misocp_minlp_first} and Table \ref{tab:compare_misocp_minlp_second}, where `Refl.' denotes the reflectance in that wavelength.

\begin{table}[H]
\centering
\caption{Comparison of MIQCP and MISOCP-based methods for Molybdenum and Niobium.}
\label{tab:compare_misocp_minlp_first}
\begin{tabular}{|c|rrrr|rrrr|}
\hline
$N=6$ & \multicolumn{4}{c|}{Molybdenum} & \multicolumn{4}{c|}{Niobium} \\ \hline
\multicolumn{1}{|l|}{} & \multicolumn{2}{c|}{MIQCP} & \multicolumn{2}{c|}{MISOCP} & \multicolumn{2}{c|}{MIQCP} & \multicolumn{2}{c|}{MISOCP} \\ \hline
\multicolumn{1}{|l|}{Wavelength(nm)} & Refl. & \multicolumn{1}{r|}{Time} & Refl. & Time & Refl. & \multicolumn{1}{r|}{Time} & Refl. & Time \\ \hline
370 & 0.995 & \multicolumn{1}{r|}{5.7} & 0.990 & 5.4 & 0.996 & \multicolumn{1}{r|}{5.0} & 0.965 & 5.4 \\
410 & 0.996 & \multicolumn{1}{r|}{7.1} & 0.923 & 5.3 & 0.996 & \multicolumn{1}{r|}{1.7} & 0.978 & 5.3 \\
450 & 0.996 & \multicolumn{1}{r|}{6.0} & 0.958 & 5.8 & 0.995 & \multicolumn{1}{r|}{5.6} & 0.907 & 5.3 \\
490 & 0.995 & \multicolumn{1}{r|}{5.5} & 0.964 & 5.5 & 0.995 & \multicolumn{1}{r|}{6.0} & 0.985 & 5.3 \\
530 & 0.995 & \multicolumn{1}{r|}{180.4} & 0.980 & 5.7 & 0.994 & \multicolumn{1}{r|}{210.3} & 0.967 & 5.5 \\
570 & 0.994 & \multicolumn{1}{r|}{92.6} & 0.963 & 5.3 & 0.994 & \multicolumn{1}{r|}{107.0} & 0.745 & 5.3 \\
610 & 0.993 & \multicolumn{1}{r|}{66.2} & 0.943 & 5.4 & 0.994 & \multicolumn{1}{r|}{98.5} & 0.950 & 5.1 \\
650 & 0.993 & \multicolumn{1}{r|}{55.0} & 0.910 & 15.2 & 0.993 & \multicolumn{1}{r|}{65.5} & 0.960 & 14.5 \\
690 & 0.993 & \multicolumn{1}{r|}{58.2} & 0.973 & 15.1 & 0.993 & \multicolumn{1}{r|}{66.0} & 0.950 & 14.7 \\
730 & 0.992 & \multicolumn{1}{r|}{64.9} & 0.971 & 15.0 & 0.993 & \multicolumn{1}{r|}{81.0} & 0.943 & 14.5 \\
770 & 0.992 & \multicolumn{1}{r|}{61.9} & 0.954 & 14.4 & 0.994 & \multicolumn{1}{r|}{87.1} & 0.987 & 14.7 \\ \hline
Average & 0.994 & \multicolumn{1}{r|}{54.9} & 0.957 & 8.9 & 0.994 & \multicolumn{1}{r|}{66.7} & 0.940 & 8.7 \\ \hline
\end{tabular}
\end{table}

\begin{table}[H]
\centering
\caption{Comparison of  MIQCP and MISOCP-based methods for Tantalum and Tungsten.}
\label{tab:compare_misocp_minlp_second}
\begin{tabular}{|c|rrrr|rrrr|}
\hline
$N=6$ & \multicolumn{4}{c|}{Tantalum} & \multicolumn{4}{c|}{Tungsten} \\ \hline
\multicolumn{1}{|l|}{} & \multicolumn{2}{c|}{MIQCP} & \multicolumn{2}{c|}{MISOCP} & \multicolumn{2}{c|}{MIQCP} & \multicolumn{2}{c|}{MISOCP} \\ \hline
Wavelength(nm) & Refl. & \multicolumn{1}{r|}{Time} & Refl. & Time & Refl. & \multicolumn{1}{r|}{Time} & Refl. & Time \\ \hline
370 & 0.995 & \multicolumn{1}{r|}{5.9} & 0.966 & 5.5 & 0.996 & \multicolumn{1}{r|}{6.1} & 0.891 & 5.5 \\
410 & 0.996 & \multicolumn{1}{r|}{5.5} & 0.973 & 5.6 & 0.996 & \multicolumn{1}{r|}{5.7} & 0.973 & 5.5 \\
450 & 0.993 & \multicolumn{1}{r|}{160.8} & 0.985 & 5.6 & 0.995 & \multicolumn{1}{r|}{151.8} & 0.977 & 5.6 \\
490 & 0.992 & \multicolumn{1}{r|}{146.1} & 0.930 & 5.4 & 0.994 & \multicolumn{1}{r|}{170.4} & 0.967 & 5.3 \\
530 & 0.991 & \multicolumn{1}{r|}{155.8} & 0.986 & 5.5 & 0.993 & \multicolumn{1}{r|}{178.6} & 0.969 & 5.5 \\
570 & 0.989 & \multicolumn{1}{r|}{134.9} & 0.844 & 5.5 & 0.993 & \multicolumn{1}{r|}{90.2} & 0.964 & 5.3 \\
610 & 0.989 & \multicolumn{1}{r|}{98.7} & 0.754 & 5.4 & 0.992 & \multicolumn{1}{r|}{65.2} & 0.867 & 5.9 \\
650 & 0.991 & \multicolumn{1}{r|}{83.8} & 0.797 & 15.1 & 0.992 & \multicolumn{1}{r|}{54.6} & 0.904 & 15.0 \\
690 & 0.993 & \multicolumn{1}{r|}{77.7} & 0.986 & 14.8 & 0.992 & \multicolumn{1}{r|}{60.8} & 0.912 & 15.0 \\
730 & 0.994 & \multicolumn{1}{r|}{73.8} & 0.964 & 14.7 & 0.991 & \multicolumn{1}{r|}{59.4} & 0.927 & 15.1 \\
770 & 0.995 & \multicolumn{1}{r|}{5.1} & 0.978 & 14.9 & 0.990 & \multicolumn{1}{r|}{67.7} & 0.944 & 15.1 \\ \hline
Average & 0.993 & \multicolumn{1}{r|}{86.2} & 0.924 & 8.9 & 0.993 & \multicolumn{1}{r|}{82.8} & 0.936 & 9.0 \\ \hline
\end{tabular}
\end{table}

Examining Table \ref{tab:compare_misocp_minlp_first} and Table \ref{tab:compare_misocp_minlp_second}, we see that by solving the MIQCP, even 6 layers of coating is sufficient to design thin films with 99\% performance when we are interested in only a single wavelength. Also, solving the MISOCP is around 6 to 10 times faster than solving the MIQCP. While achieving faster solution times, the MISOCP-based method provides designs with 5 to 8\% worse performance compared to the MIQCP-based method.

\subsection{Experiments for the Visible Spectrum}\label{subsec:experiment_vis}


We solve Problem~\eqref{eq:discrete-nonconvex} with its objective function discretized at different intervals and report the results in Table~\ref{tab:vis_spectrum_results}, where
$\Lambda(10):= \{370,380,\dots,770\}$, $\Lambda(20):= \{370,390,\dots,770\}$ and $\Lambda(40):= \{370,410,\dots,770\}$. The allowed thicknesses for the coatings are the same as in Section~\ref{sec:singlewavelength}. Then, we report the average reflectance over the visible spectrum in Table~\ref{tab:vis_spectrum_results} and the solution times of the MIQCP problems in Table~\ref{tab:vis_spectrum_results_time}, where `TL' denotes that the solver reached the time limit of 5 hours.

\begin{table}[H]
\centering
\caption{Average reflectance rates for the visible spectrum for the MIQCP-based method.}
\label{tab:vis_spectrum_results}
\begin{tabular}{|c|ccc|ccc|ccc|}
\hline
\multicolumn{1}{|l|}{\textbf{}} & \multicolumn{3}{c|}{\textbf{$ \Lambda(10)$}} & \multicolumn{3}{c|}{\textbf{$\Lambda(20)$}} & \multicolumn{3}{c|}{\textbf{$\Lambda (40)$}} \\ \hline
\textbf{Substrate/$N$} & \textbf{6} & \textbf{10} & \textbf{14} & \textbf{6} & \textbf{10} & \textbf{14} & \textbf{6} & \textbf{10} & \textbf{14} \\ \hline
\textbf{Molybdenum} & 0.939 & 0.972 & \textbf{0.993}& 0.940 & 0.970 & 0.992 & 0.928 & 0.960 & 0.984 \\
\textbf{Niobium} & 0.942 & 0.979 & \textbf{0.994}& 0.941 & 0.971 & 0.991 & 0.937 & 0.948 & 0.985 \\
\textbf{Tantalum} & 0.931 & 0.976 & \textbf{0.991}& 0.927 & 0.977 & 0.991 & 0.916 & 0.961 & 0.979 \\
\textbf{Tungsten} & 0.924 & 0.966 & \textbf{0.993}& 0.924 & 0.954 & 0.982 & 0.905 & 0.954 & 0.984 \\ \hline
\end{tabular}
\end{table}

\begin{table}[H]
\centering
\caption{Runtimes for the MIQCP-based method (seconds).}
\label{tab:vis_spectrum_results_time}
\begin{tabular}{|c|rrr|rrr|rrr|}
\hline
\multicolumn{1}{|l|}{\textbf{}} & \multicolumn{3}{c|}{\textbf{$ \Lambda(10)$}} & \multicolumn{3}{c|}{\textbf{$\Lambda(20)$}} & \multicolumn{3}{c|}{\textbf{$\Lambda (40)$}} \\ \hline
\textbf{Substrate/$N$} & \multicolumn{1}{c}{\textbf{6}} & \multicolumn{1}{c}{\textbf{10}} & \multicolumn{1}{c|}{\textbf{14}} & \multicolumn{1}{c}{\textbf{6}} & \multicolumn{1}{c}{\textbf{10}} & \multicolumn{1}{c|}{\textbf{14}} & \multicolumn{1}{c}{\textbf{6}} & \multicolumn{1}{c}{\textbf{10}} & \multicolumn{1}{c|}{\textbf{14}} \\ \hline
\textbf{Molybdenum} & TL & TL & TL & TL & TL & 524 & TL & TL & 223 \\
\textbf{Niobium} & TL & TL & TL & TL & TL & TL & TL & TL & 2701 \\
\textbf{Tantalum} & TL & TL & TL & TL & TL & TL & TL & TL & 483 \\
\textbf{Tungsten} & TL & TL & TL & TL & TL & TL & TL & TL & TL \\ \hline
\end{tabular}
\end{table}

\begin{figure}[H]
\centering
    \begin{subfigure}{0.45\linewidth}
\includegraphics[width=\linewidth]{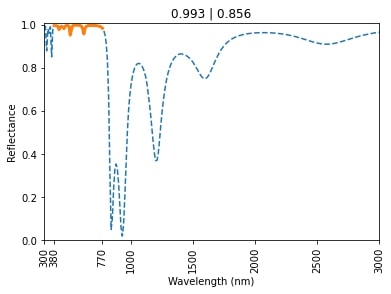}
\caption{Molybdenum}
\label{fig:subfig1}
    \end{subfigure}\hfil
    \begin{subfigure}{0.45\linewidth}
\includegraphics[width=\linewidth]{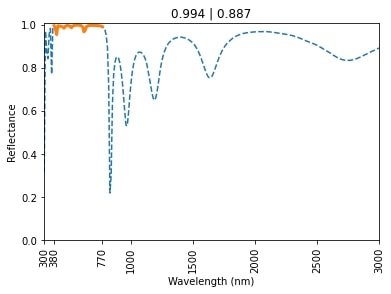}
\caption{Niobium}
\label{fig:subfig2}
    \end{subfigure}

\smallskip
     \begin{subfigure}{0.45\linewidth}
\includegraphics[width=\linewidth]{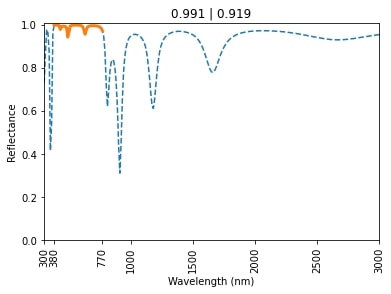}
\caption{Tantalum}
\label{fig:subfig3}
    \end{subfigure}\hfil
     \begin{subfigure}{0.45\linewidth}
\includegraphics[width=\linewidth]{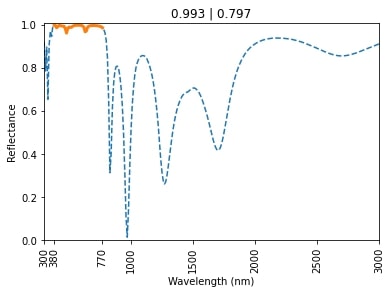}
\caption{Tungsten}
\label{fig:subfig4}
    \end{subfigure}
    \caption{Reflectance profiles with $N=14$ layers optimized for the visible spectrum for the MIQCP-based method.}
\label{fig:vis_reflectance_profiles}
\end{figure}

Examining Table~\ref{tab:vis_spectrum_results}, we see that the best average reflectance rates for the visible spectrum are achieved when $N=14$ layers are coated under the objective function discretization defined by $\Lambda(10)$, the densest setting. Since the natural upper bound on the reflectance rate is 1, the results are quite satisfactory.
On one hand, it is expected for the setting with $\Lambda(10)$ to perform the best  since the discretization in the other settings are coarser, the optimization model does not consider the reflectance at some wavelengths. On the other hand, coarser discretization may lead to faster convergence for larger $N$ values. For example, for $\Lambda(40)$, $N=14$ and the substrates Tantalum, Niobium and Molybdenum, the optimization model converges within less than 45 minutes but for all other problem instances (except Molybdenum at $N=14$ and $\Lambda(20)$, which takes around 9 minutes to solve), the solver terminates upon reaching the time limit. We present the reflectance profiles of the best solutions found by the solver in Figure~\ref{fig:vis_reflectance_profiles}. 

The results suggest that while 14 layers are enough to get the average reflectance rate beyond 99\% for the visible spectrum, the broad spectrum performances are not satisfactory. This outcome is expected as the objective function is limited to the visible spectrum. To improve reflectance performance across a broader range, it is necessary to extend the objective function to include wavelengths beyond the visible region.

We also carry out the same experiments by utilizing the MISOCP-based method. After solving the instances with MISOCPs with a 30-minute time limit, we evaluate their reflectance value as in Table~\ref{tab:vis_spectrum_MISOCP_results} and report the time it takes to solve the MISOCP in Table~\ref{tab:vis_spectrum_MISOCP_time}. The results show that solving the MISOCPs can be 10 times faster compared to solving the MIQCPs. However, solving the MISOCPs yield on average 10.5\% worse average reflectance over all instances, compared to the solutions obtained by solving the MIQCPs. The speed of the MISOCP is promising as number of coating layers increase. Therefore, it might be helpful for the broad spectrum, where we need more coating materials to achieve good reflectance values as demonstrated in Section~\ref{subsec:experiment_broad}.

\begin{table}[H]
\centering
\caption{Average reflectance rates for the visible spectrum for the MISOCP-based method.}
\label{tab:vis_spectrum_MISOCP_results}
\begin{tabular}{|c|ccc|ccc|ccc|}
\hline
\multicolumn{1}{|l|}{\textbf{}} & \multicolumn{3}{c|}{\textbf{$ \Lambda(10)$}} & \multicolumn{3}{c|}{\textbf{$\Lambda(20)$}} & \multicolumn{3}{c|}{\textbf{$\Lambda (40)$}} \\ \hline
\textbf{Substrate/$N$} & \textbf{6} & \textbf{10} & \textbf{14} & \textbf{6} & \textbf{10} & \textbf{14} & \textbf{6} & \textbf{10} & \textbf{14} \\ \hline
\textbf{Molybdenum} & 0.796 & 0.892 & 0.915 & 0.783 & 0.847 & 0.911 & 0.814 & 0.884 & 0.864 \\
\textbf{Niobium} & 0.827 & 0.931 & 0.887 & 0.840 & 0.866 & 0.883 & 0.874 & 0.843 & 0.890 \\
\textbf{Tantalum} & 0.775 & 0.803 & 0.885 & 0.784 & 0.888 & 0.918 & 0.804 & 0.852 & 0.916 \\
\textbf{Tungsten} & 0.790 & 0.824 & 0.907 & 0.784 & 0.874 & 0.875 & 0.823 & 0.868 & 0.932 \\ \hline
\end{tabular}
\end{table}

\begin{table}[H]
\centering
\caption{Runtimes for the MISOCP-based method (seconds).}
\label{tab:vis_spectrum_MISOCP_time}
\begin{tabular}{|c|rrr|rrr|rrr|}
\hline
\multicolumn{1}{|l|}{\textbf{}} & \multicolumn{3}{c|}{\textbf{$ \Lambda(10)$}} & \multicolumn{3}{c|}{\textbf{$\Lambda(20)$}} & \multicolumn{3}{c|}{\textbf{$\Lambda (40)$}} \\ \hline
\textbf{Substrate/$N$} & \multicolumn{1}{c}{\textbf{6}} & \multicolumn{1}{c}{\textbf{10}} & \multicolumn{1}{c|}{\textbf{14}} & \multicolumn{1}{c}{\textbf{6}} & \multicolumn{1}{c}{\textbf{10}} & \multicolumn{1}{c|}{\textbf{14}} & \multicolumn{1}{c}{\textbf{6}} & \multicolumn{1}{c}{\textbf{10}} & \multicolumn{1}{c|}{\textbf{14}} \\ \hline
\textbf{Molybdenum} & 1800 & 1721 & 1800 & 1800 & 299 & 703 & 1800 & 98 & 125 \\
\textbf{Niobium} & 1800 & 1483 & 1800 & 1800 & 403 & 325 & 1293 & 72 & 148 \\
\textbf{Tantalum} & 1800 & 1244 & 1800 & 169 & 419 & 518 & 1800 & 72 & 85 \\
\textbf{Tungsten} & 1800 & 1800 & 1800 & 237 & 175 & 624 & 1800 & 56 & 75 \\ \hline
\end{tabular}
\end{table}

\subsection{Experiments for the Broad Spectrum}\label{subsec:experiment_broad}

To enhance the reflectance for the broad spectrum (300-3000 nm), we solve Problem~\eqref{eq:discrete-nonconvex} with different settings. We experiment with different discretizations of the objective function and material thicknesses. As observed in Figure~\ref{fig:vis_reflectance_profiles}, the solutions obtained by focusing on the visible spectrum exhibit high reflectance for wavelengths beyond 1500 nm, even though the optimization targets the visible range. Motivated by this observation, the setting $\Lambda_2$ considers wavelengths in the 300–1500 nm range.

In addition, we investigate two different sets of admissible material thicknesses. In the first set, thickness values are close to the quarter wavelength optical thicknesses corresponding to the visible spectrum. In the second set, thicknesses are close to the quarter wavelength optical thicknesses over the broad spectrum. The settings are specified in Tables~\ref{tab:broad_discretizations} and \ref{tab:broad_thicknesses}.

\begin{table}[H]
\centering
\caption{Settings for objective function discretization.}
\label{tab:broad_discretizations}
\begin{tabular}{|c|l|}
\hline
\textbf{Setting} & \textbf{Wavelengths of Light (nm)} \\ \hline
$\Lambda_1$&
$\{300,340,\dots,1500\} \cup \{1750,2000,\dots,3000\}$ \\ \hline
$\Lambda_2$&
$\{300,340,\dots,1460,1500\}$ \\ \hline
\end{tabular}
\end{table}

\begin{table}[H]
\centering
\caption{Settings for material thicknesses.}
\label{tab:broad_thicknesses}
\begin{tabular}{|c|cc|}
\hline
\textbf{Setting} & \textbf{Coating Material} & \textbf{Material Thicknesses (nm)} \\ \hline
\multirow{2}{*}{\textbf{$\Theta^1$}} & \ch{MgF2} & $\{50,60,\dots,280\}$ \\
 & \ch{TiO2} & $\{20,30,\dots,140\}$ \\ \hline
\multirow{2}{*}{\textbf{$\Theta^2$}} & \ch{MgF2} & $\{50,70,\dots,550\}$ \\
 & \ch{TiO2} & $\{20,40,\dots,300\}$ \\ \hline
\end{tabular}
\end{table}

Having defined the experimental setup, we are now able to present the computational results of the MIQCP-based method in Table~\ref{tab:broad_spectrum_results}, where `N/A' denotes that the solver is not able to find a feasible solution within the time limit.

\begin{table}[H]
\centering
\caption{Average reflectance rates for the broad spectrum for the MIQCP-based method with $N=20$ layers of coating.}
\label{tab:broad_spectrum_results}
\begin{tabular}{|c|c|rr|rr|}
\hline
\multicolumn{1}{|l|}{} & \multicolumn{1}{l|}{} & \multicolumn{2}{c|}{\textbf{1-hour time limit}} & \multicolumn{2}{c|}{\textbf{5-hour time limit}} \\ \hline
\textbf{Substrate} & \textbf{Setting} & \textbf{$\Lambda_1$} & \textbf{$\Lambda_2$} & \textbf{$\Lambda_1$} & \textbf{$\Lambda_2$} \\ \hline
\multirow{2}{*}{\textbf{Molybdenum}} & \textbf{$\Theta^1$} & 0.915 & 0.941 & 0.942 & 0.951 \\
 & \textbf{$\Theta^2$} & 0.947 & 0.945 & \textbf{0.964} & 0.962 \\ \hline
\multirow{2}{*}{\textbf{Niobium}} & \textbf{$\Theta^1$} & N/A & 0.946 & 0.939 & 0.947 \\
 & \textbf{$\Theta^2$} & \textbf{0.958} & 0.950 & 0.954 & 0.950 \\ \hline
\multirow{2}{*}{\textbf{Tantalum}} & \textbf{$\Theta^1$} & 0.909 & 0.953 & \textbf{0.967} & 0.961 \\
 & \textbf{$\Theta^2$} & 0.955 & 0.954 & \textbf{0.966} & \textbf{0.966} \\ \hline
\multirow{2}{*}{\textbf{Tungsten}} & \textbf{$\Theta^1$} & N/A & 0.909 & 0.943 & 0.949 \\
 & \textbf{$\Theta^2$} & 0.917 & 0.942 & \textbf{0.957} & 0.939 \\ \hline
\end{tabular}
\end{table}

\begin{figure}[H]
\centering
    \begin{subfigure}{0.45\linewidth}
\includegraphics[width=\linewidth]{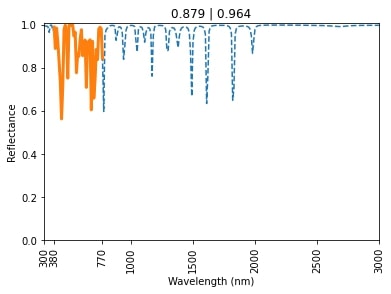}
\caption{Molybdenum}
\label{fig:subfig1b}
    \end{subfigure}\hfil
    \begin{subfigure}{0.45\linewidth}
\includegraphics[width=\linewidth]{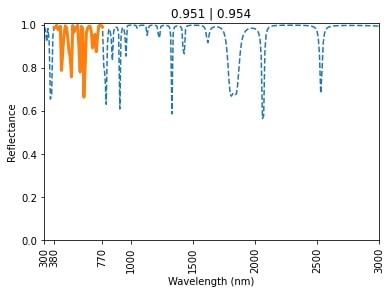}
\caption{Niobium}
\label{fig:subfig2b}
    \end{subfigure}

\smallskip
     \begin{subfigure}{0.45\linewidth}
\includegraphics[width=\linewidth]{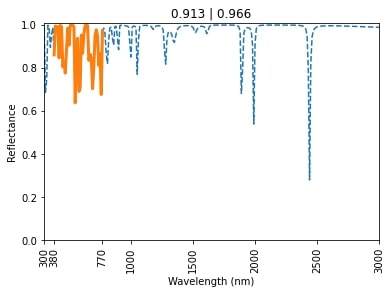}
\caption{Tantalum}
\label{fig:subfig3b}
    \end{subfigure}\hfil
     \begin{subfigure}{0.45\linewidth}
\includegraphics[width=\linewidth]{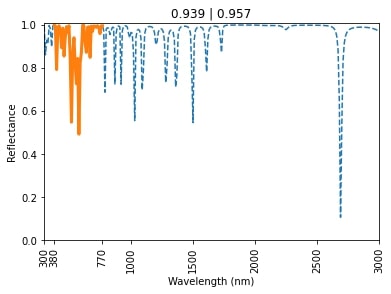}
\caption{Tungsten}
\label{fig:subfig4b}
    \end{subfigure}
    \caption{Reflectance profiles with $N=20$ layers optimized for the broad spectrum.}
\label{fig:broad_reflectance_profiles}
\end{figure}

Examining the results of Table~\ref{tab:broad_spectrum_results}, we see that we are able to achieve at least 95.8\% reflectance in the broad spectrum with $N=20$ layers of coating. While the average reflectance can be further improved by increasing the number of layers, doing so makes the problem more challenging to solve. The results indicate that the  setting $\Lambda_1$ with the set of material thicknesses defined by $\Theta^2$ generally outperforms the other combinations. Therefore, we present the reflectance profiles of the setting in Figure~\ref{fig:broad_reflectance_profiles}.

Examining Figure~\ref{fig:broad_reflectance_profiles}, we observe that when the coating design is optimized for the broad spectrum, we may see a deterioration in the reflectance for the visible spectrum. Since the broad spectrum is much wider compared to the visible spectrum and the coating has only 20 layers, this deterioration is expected. For the Molybdenum substrate, while we are able to achieve 96.4\% performance for the broad spectrum, the visible spectrum performance drops below 88\%. While this trade off is natural, further improvements in reflectance across both spectra would require increasing the number of layers.

To assess the success of the proposed MISOCP-based method, we solve the MISOCP problems with a time limit of one hour. The results of the experiments are presented in Table~\ref{tab:broad_spectrum_results_MISOCP}. 

\begin{table}[H]
\centering
\caption{Average reflectance rates for the broad spectrum for the MISOCP-based method with $N=20$ layers of coating.}
\label{tab:broad_spectrum_results_MISOCP}
\begin{tabular}{|c|c|rr|}
\hline
\textbf{Substrate} & \textbf{Setting} & \textbf{$\Lambda_1$} & \textbf{$\Lambda_2$} \\ \hline
\multirow{2}{*}{\textbf{Molybdenum}} & \textbf{$\Theta^1$} & N/A & \textbf{0.929} \\
 & \textbf{$\Theta^2$} & 0.893 & 0.919 \\ \hline
\multirow{2}{*}{\textbf{Niobium}} & \textbf{$\Theta^1$} & N/A & 0.938 \\
 & \textbf{$\Theta^2$} & 0.912 & \textbf{0.952} \\ \hline
\multirow{2}{*}{\textbf{Tantalum}} & \textbf{$\Theta^1$} & 0.957 & 0.964 \\
 & \textbf{$\Theta^2$} & N/A &\textbf{ 0.966} \\ \hline
\multirow{2}{*}{\textbf{Tungsten}} & \textbf{$\Theta^1$} & 0.910 & \textbf{0.945} \\
 & \textbf{$\Theta^2$} & N/A & 0.943 \\ \hline
\end{tabular}
\end{table}

Comparing the reflectance values reported in  Table~\ref{tab:broad_spectrum_results}  and Table~\ref{tab:broad_spectrum_results_MISOCP}, we observe that the MISOCP-based method produces high-quality solutions. In particular, for Tantalum and Niobium substrates, the MISOCP-based approach yields solutions of nearly identical quality within one hour. The ability of the MISOCP-based method to produce high-quality solutions within reasonable computational times highlights its potential as a promising approach.

In order to demonstrate the effectiveness of the proposed methods, we compare our results with a heuristic method from the literature \citep{keccebacs2018enhancing}. In this paper, the authors first specify nine wavelengths: 450, 500, 750, 900, 1000, 1200, 1500, 2000, 2200 nm. Then, for each wavelength, they design 7-layer films using the optical quarter wavelength thicknesses of the materials and then stack them on top of the substrate. Since this approach results in films with 63 layers, to obtain comparable solutions with our study, we utilize their heuristic with different number of layers. We report the average reflectance values in the broad spectrum for different number of layers in Table~\ref{tab:kecebascomparison}. For instance, KŞ-9$\times$2 means that 2-layer films are designed for each wavelength and then stacked, resulting in 18 layers.

\begin{table}[H]
\centering
\caption{Reflectance values of proposed methods and the heuristic.}
\label{tab:kecebascomparison}
\begin{tabular}{|c|c|c|c|c|}
\hline
\textbf{Method-$N$} & \textbf{Molybdenum} & \textbf{Niobium} & \textbf{Tantalum} & \textbf{Tungsten} \\ \hline
\textbf{MIQCP-20} & 0.964 & 0.954 & 0.967 & 0.957 \\
\textbf{MISOCP-20} & 0.929 & 0.952 & 0.966 & 0.945 \\ \hline
\textbf{KŞ-9$\times$2} & 0.943 & 0.943 & 0.955 & 0.924 \\
\textbf{KŞ-9$\times$3} & 0.941 & 0.943 & 0.953 & 0.925 \\
\textbf{KŞ-9$\times$4} & 0.979 & 0.974 & 0.980 & 0.970 \\
\textbf{KŞ-9$\times$5} & 0.978 & 0.976 & 0.978 & 0.973 \\
\textbf{KŞ-9$\times$6} & 0.991 & 0.986 & 0.992 & 0.986 \\
\textbf{KŞ-9$\times$7} & 0.988 & 0.987 & 0.990 & 0.985 \\ \hline
\end{tabular}
\end{table}

Examining Table~\ref{tab:kecebascomparison}, we observe that solving the MIQCP-based optimization model for 20 layers achieves much better reflectance values than the heuristic with 27 layers. The heuristic method surpasses the optimization based method when we allow for 36 layers to be coated, which is natural because 16 more layers are allowed to enhance the reflectance. On the other hand, MISOCP-based model also surpasses the heuristic method with 27 layers except for Molybdenum substrate, which is promising.

\section{Conclusion}\label{sec:conclusion}

In this work, we develop an MIQCP-based global optimization formulation and an MISOCP-based relaxation for the Multi-Layer Thin Films Problem. We design coatings that remain effective even when the wavelength of the incident light varies. Our extensive computational experiments with the MIQCP-based method yield solutions that provide 99\% reflectance in the visible spectrum and 95\% reflectance in the broad spectrum. Our proposed MISOCP-based method yields competitive results with both the MIQCP-based model and a heuristic method from the literature, and is promising to yield good quality solutions as the number of coating layers increase.

There are two potential research avenues to explore: Firstly, further improvements for broad spectrum reflectance could be obtained by considering designs with a larger number of coating layers, provided that reasonable computational times can be maintained. Secondly, the tightness of the MISOCP relaxation can be enhanced and   its scalability can be improved as the number of layers increases.

\section*{Disclosure statement}

No potential conflict of interest was reported by the author(s).

\section*{Funding}
This work was supported by the  Scientific and Technological Research Council of Turkey with grant number 120C151.

\section*{Data availability statement}
The data relevant to calculation of transfer matrices are obtained from \cite{keccebacs2018enhancing}, taken from multiple sources \citep{malitson1965interspecimen,palik1998handbook,dodge1984refractive,dodge2refractive,golovashkin1969optical}.

\section*{Author contributions}
CRediT: \textbf{Deniz Tuncer:} Conceptualization, Data curation, Formal analysis, Investigation, Methodology, Project administration, Software, Visualization, Writing - original draft;
\textbf{Burak Kocuk:} Conceptualization, Funding acquisition, Methodology, Resources, Software, Supervision, Writing - review \& editing.

\section*{ORCID}

\textit{Deniz Tuncer} https://orcid.org/0000-0003-0176-4937

\noindent\textit{Burak Kocuk} https://orcid.org/0000-0002-4218-1116

\bibliographystyle{tfcad}
\bibliography{references}

\end{document}